# Brief Research Notes on Transformation of Series and Special Functions

**Nikos Bagis**
Department of Informatics Aristotle University
Thessaloniki, Greece
bagkis@hotmail.com

**Abstract**
In this work we derive results concerning Elliptic Functions using as tools general formulas from previous work.

## 1. Introduction

**Theorem 1. (see [G,B]).** Let $f$ be analytic function in the upper half plane $\mathrm{Im}(z)>0$ and continuous in $\mathrm{Im}(z)\geq 0$ and there exists $C, N>0$ and $0\leq\beta<\pi$ such that

$$|f(z)|\leq C(1+|z|)^N e^{\beta|\mathrm{Re}(z)|}, \text{ for every } z \text{ in } \mathrm{Im}(z)\geq 0$$

Let also $a>0$ and $f$ is odd, then

$$\sqrt{a}\left(\frac{f'(0)}{2\pi}+\sum_{n=1}^{\infty}\frac{f(na)}{\sinh(\pi na)}\right)=\sqrt{\frac{b2}{\pi}}\left(c_o/2+i\sum_{n=1}^{\infty}\frac{(-1)^n f(ni)}{e^{bn}-1}\right)$$

when $ab=2\pi$ and $c_o=\lim_{x\to 0^+}\left(\sum_{n=1}^{\infty}(-1)^n i f(ni)e^{-nx}\right)$.

**Theorem 2. (submitted on a journal)**
If $a>0$, $b>0$, $ab=2\pi$, then

$$a\lim_{r\to 1^{-1}}\left(\sum_{n=1}^{\infty}y(n)f(n)r^n\right)+a\sum_{k=1}^{\infty}\frac{y(n)f(n)}{e^{an}-1}+$$

$$+a\sum_{n=1}^{\infty}\left(\sum_{m=1}^{\infty}f(m)y(m)e^{-mna}\right)-c-2\sum_{n=1}^{\infty}\varphi(nb)=0$$

Where $\varphi(x)=\mathrm{Re}\left(M\Psi(ix)f(-ix)\right)$, $c=\lim_{h\to 0}\mathrm{Re}\left(M\Psi(ih)f(-ih)\right)$, $\Psi(x)=\sum_{n=0}^{\infty}y(n)x^n$ and $f,\Psi,x$, real, $f(0)=0$. Provided that all sums converge.

## 2. Propositions

**Proposition 1.**

$$\sum_{n=1}^{\infty}\frac{\cosh(2tn)}{n\sinh(\pi an)}=\log(P_0)-\log(\vartheta_4(it,e^{-a\pi})) \qquad :(1)$$

Where $P_0=\prod_{n=1}^{\infty}\left(1-e^{-2na\pi}\right)$ and $\vartheta_4(u,q)=1+2\sum_{n=1}^{\infty}(-1)^n q^{n^2}\cos(2nu)$.

**Proof.** From [An] pg.170 relation (13-2-12) and the definition of Theta functions we get

$$\vartheta_4(z,q)=P_0(q)\prod_{n=0}^{\infty}(1-q^{2n-1}e^{2iz})(1-q^{2n-1}e^{-2iz})$$

Where

$$P_0(q)=\prod_{n=1}^{\infty}(1-q^{2n}),\ P_0:=\prod_{n=1}^{\infty}(1-e^{-2n\pi a})$$

By taking the logarithm of both sides and expanding the logarithm of the individual terms in a power series, it is simple to show that

$$\sum_{n=1}^{\infty}\frac{\cosh\left(2nt\right)}{n\sinh(\pi an)}=\log\left(P_0\right)-\log\left(\vartheta_4(it,q)\right)$$

Where $q=e^{-\pi a}$, $a>0$. □

**Proposition 2.**

$$4\sum_{n=1}^{\infty}\frac{(-1)^n\sin(\theta n)^2}{(e^{2\pi na}-1)n}=\log\left(\frac{\vartheta_4(i\theta/a,e^{-\pi/a})}{\vartheta_4(0,e^{-\pi/a})\cos(\theta)}\right)-\frac{\theta^2}{a\pi} \qquad :(2)$$

where

**Proof.** From Theorem 1 with $f(t)=\frac{\cos(2tc)-1}{t}$, we get

$$-ac^2/\pi+\log(\cosh(c))-4\sum_{n=1}^{\infty}\frac{(-1)^n\sinh^2(cn)}{(e^{2\pi n/a}-1)n}=\sum_{n=1}^{\infty}\frac{1}{n\sinh(an\pi)}-\sum_{k=1}^{\infty}\frac{\cos(2acn)}{k\sinh(an\pi)}$$

using relation (1) we get the result. □

Differentiating (7) with respect to $c$ and letting $c\to\pi$ we easily get the next formula:

$$\Theta(iz, e^{-z}) + 2i\vartheta_4(iz, e^{-z}) = 0$$

where $\Theta(u,z) = \dfrac{\partial \vartheta_4(u,z)}{\partial u}$

**Proposition 3.**

$$2\frac{\partial^2}{\partial t^2}\log\left(\vartheta_4\left(\frac{it\pi}{2}, e^{-2\pi a}\right)\right)\Bigg|_{t=a} = K(k_a)E(k_a) - K(k_a)^2 \qquad :(3)$$

Where $k_a$ is the solution of the equation $\dfrac{K\left(\sqrt{1-k_a^2}\right)}{K(k_a)} = a$,

$$K(x) = \int_0^{\pi/2} \frac{1}{\sqrt{1-x^2\sin^2(\theta)}}d\theta$$

is the elliptic integral of the first kind and $E(x) = \int_0^{\pi/2} \sqrt{1-x^2\sin^2(\theta)}d\theta$ $a > 0$ is the elliptic integral of the second kind.

**Proof.** The following relations hold

$$\sum_{n=1}^{\infty} \frac{(-1)^n n}{e^{2\pi n/a}-1} = \frac{1}{8} - \frac{a}{4\pi} + \frac{a^2 K(k_a)}{2\pi^2}(E(k_a) - K(k_a)) \qquad :(4)$$

(see [G,B])

From Theorem 2, setting $y(n) = (-1)^n$, $f(n) = n$ and observing that

$$\sum_{n=1}^{\infty} \frac{1}{\cosh^2(\pi n x)} = -4\sum_{n=1}^{\infty} \frac{(-1)^n n}{e^{2\pi n x}-1}$$

We have

$$-\frac{1}{4} + \frac{a}{2\pi} + 2\sum_{n=1}^{\infty} \frac{(-1)^n n}{e^{2n\pi/a}-1} + 2a^2 \sum_{n=1}^{\infty} \frac{n\cosh(an\pi)}{\sinh(2an\pi)} = 0 \qquad :(5)$$

Combining (4), (1), (5) we get

$$2\frac{\partial^2}{\partial t^2}\log\left(\vartheta_4\left(\frac{it\pi}{2}, e^{-2\pi a}\right)\right)\Bigg|_{t=a} = K(k_a)E(k_a) - K(k_a)^2,$$

whenever $\dfrac{K'(k_a)}{K(k_a)} = a$. □

**Proposition 4.** If $a > 0$ and $x$ real then

$$\partial_x \left( e^{\frac{x^2 a}{\pi}} \frac{\vartheta_2(x, e^{-\pi/a})}{\vartheta_4(iax, e^{-a\pi})} \right) = 0 \qquad :(6)$$

**Proof.** Let $y(n)=\frac{(-1)^n\sin(nv)}{n}$ and $f(n)=n$ in Theorem 2 then we have $M\Psi(s)=\frac{\pi\csc(\pi s)\sin(\pi s)}{2s}$. Hence if $a,v>0$, noting that

$$\sum_{n=1}^{\infty}\frac{(-1)^n\sin(nv)}{e^{an}-1}=-1/2\sum_{n=1}^{\infty}\frac{\sin(v)}{\cos(v)+\cosh(an)}$$

we have

$$\frac{a}{2}\tan\left(\frac{v}{2}\right)=v+2a\sum_{n=1}^{\infty}\frac{(-1)^n\sin(nv)}{e^{an}-1}+2\pi\sum_{n=1}^{\infty}\operatorname{csch}\left(\frac{2n\pi^2}{a}\right)\sinh\left(\frac{2\pi nv}{a}\right)\qquad :(7)$$

It is known that

$$4\sum_{n=1}^{\infty}\frac{(-1)^n\sin(2nz)q^{2n}}{1-q^{2n}}=\tan(z)+\frac{1}{\vartheta_2(z,q)}\frac{\partial\vartheta_2(z,q)}{\partial z}\qquad :(8)$$

where $\vartheta_2(z,q)=\sum_{n=-\infty}^{\infty}q^{(n+1/2)^2}e^{(2n+1)iz}$ is the known theta function.

Thus

$$2\pi\sum_{n=1}^{\infty}\frac{\sinh\left(2n\pi za\right)}{\sinh\left(n\pi^2a\right)}=-2z-\frac{1}{a}\partial_z\log(\vartheta_2(z,e^{-1/a}))$$

If we differentiate (1) we get the relation between the two theta functions. □

**Proposition 5.**

$$-2\sum_{n=1}^{\infty}\frac{e^{2n\pi/b}}{(1+e^{2n\pi/b})^3}+\sum_{n=1}^{\infty}\frac{(-1)^n n^2}{e^{2n\pi/b}-1}=\frac{1}{8}-\frac{b}{4\pi}+\frac{b^2}{2\pi^2}\left(E(k_b)K(k_b)-K(k_b)^2\right)\qquad :(9)$$

**Proof.** Form Theorem 2 we have

$$-1-8\sum_{n=1}^{\infty}\frac{e^{2n\pi/b}}{(1+e^{2n\pi/b})^3}+4\sum_{n=1}^{\infty}\frac{(-1)^n n(n+1)}{e^{2n\pi/b}-1}+\frac{2b}{\pi}=-4b^2\sum_{n=1}^{\infty}\frac{n}{\sinh(\pi nb)}$$

But when $q=e^{-\pi b}$ then holds

$$2\sum_{n=1}^{\infty}\frac{n}{q^{-n}-q^n}=\sum_{n=1}^{\infty}\frac{n}{\sinh(\pi bn)}=\frac{K(k_b)}{\pi^2}(K(k_b)-E(k_b))$$

Using (4) we get the result. □

**Proposition 6.**

If $a$ is appositive real number then

$$\partial_z\left(\vartheta_2(z,e^{-\pi/(2a)})\right)\Big|_{z=\pi/4}=-\frac{2a}{\pi}\vartheta_2\left(\pi/4,e^{-\pi/(2a)}\right)K(k_a) \qquad :(10)$$

**Proof.**

From (8) using the identity: sinh(2$x$) = 2 sinh($x$) cosh($x$) and the fact that

$\sum_{n=1}^{\infty}\frac{1}{\cosh(n\pi a)}=\frac{1}{2}+\frac{K(k_a)}{\pi}$, we get easily the result. □

**Proposition 7.** If $K(k)=\int_0^{\pi/2}\frac{1}{\sqrt{1-k^2\sin^2(\theta)}}d\theta$ and $\frac{K(1-k)}{K(k)}=a$ then

$$da=\frac{K'(k)}{E(k)K(k)-K(k)^2}dk \qquad :(11)$$

**Proof.** It holds

$$2\frac{\partial^2}{\partial t^2}\log\left(\vartheta_4\left(\frac{it\pi}{2},e^{-2\pi a}\right)\right)\Big|_{t=a}=K(k_a)E(k_a)-K(k_a)^2$$

From (1) we have $w(x)=2\frac{\partial}{\partial t}\log\left(\vartheta_4\left(\frac{it\pi}{2},e^{-2\pi x}\right)\right)\Big|_{t=x}=-\pi\sum_{n=1}^{\infty}\frac{1}{\cosh(n\pi x)}=\frac{\pi}{2}-K(k_x)$

then exists function $A(k)$ such that

$$w(r_2)-w(r_1)=-\int_{k_1}^{k_2}\left(K(k)E(k)-K(k)^2\right)A(k)dk$$

and it is known that $K'(k)=\frac{dK(k)}{dk}=\frac{E(k)-(1-k)K(k)}{2(1-k)k}$

Combining the above relations we get (11). □

**Proposition 8.**

Set $\{r,k\}=\frac{dr}{dk}=\frac{K'(k)}{E(k)K(k)-K(k)^2}$ and $q=e^{-\pi r}$ then

$$24\sum_{n=1}^{\infty}\frac{nq^n}{1-q^n}=1+\frac{6E(k_r)+(k_r-5)K(k_r)}{\pi k_r(1-k_r)K(k_r)\{r,k\}}$$

**Proof.**

From $f(-q)=\frac{2^{1/3}}{\pi^{1/2}}q^{-1/24}k^{1/24}(1-k)^{1/6}K(k)^{1/2}$ (see [W,W]) taking the logarithmic derivative with respect to *k* and using Proposition 7 we get the result. □

**Proposition 9. (observation)**

$$\frac{1}{\sqrt{1-x}}K\left(\frac{x}{x-1}\right)=K(x)$$

**Proposition 10.**

Let $F$ be even and analytic in the real line with $F(0)=F'(0)=F''(0)=0$. Then for real $a>0$, $ab=2\pi$ we have:

$$2\int_1^2\left(\frac{1}{t}\sum_{n=1}^{\infty}G\left(\frac{t}{2\pi in}\right)\right)dt+2\sum_{n=1}^{\infty}\frac{(-1)^nF(n)}{n(e^{an}-1)}-\sum_{n=1}^{\infty}\frac{F(ibn)}{n\sinh(bn\pi)}=0$$

where $F(x)=\sum_{k=1}^{\infty}f_{2n}x^{2n}$, $G(x)=\sum_{n=1}^{\infty}g_{2n}x^{2n}$, with $f_n=\frac{g_n}{n!}$

**Proof.** From Theorem 2 with $f(x)=x^{2v}$ and $y(n)=\frac{(-1)^n}{n}$ we have

$$\sum_{v=1}^{\infty}\frac{f_e^{(2v)}(0)}{(2v)!}(2^{2v}-1)\zeta(1-2v)+2\sum_{n=1}^{\infty}\frac{(-1)^nf_e(n)}{n(e^{an}-1)}-\sum_{n=1}^{\infty}\frac{f_e(ibn)}{n\sinh(bn\pi)}=0$$

when $ab=2\pi$.

Let $B_n$ denotes the Bernoulli numbers and $\zeta$ is the Riemann's Zeta function. Use the relations

$$\sum_{n=1}^{\infty}\left(f\left(\frac{x}{2\pi in}\right)+f\left(\frac{-x}{2\pi in}\right)-2f(0)\right)=-\sum_{n=1}^{\infty}\frac{f^{(2n)}(0)}{(2n)!}\frac{B_{2n}}{(2n)!}x^{2n},$$

$$\zeta(1-2v)+\frac{B_{2v}}{2v}=0$$

([Ba,Gl] and [A,S] respectively), and the result will follows. □

**Lemma.**

Let $f(x)=\sum_{n=0}^{\infty}f_nx^n$ be analytic function on $\mathbb{R}$ such that for every $a,b>0$ there exist constant $M_f$:

$$|f(t)|\le M_f\left(1+|t|^a\right)e^{-|t|b}$$

Let also $\varphi:\mathbb{R}\to\mathbb{C}$, such that for every $c>0$ there exist constant $M_\varphi$:

$$|\varphi(t)| \le M_\varphi |t|^c$$

then

$$\int_0^\infty f(t)\varphi(t)e^{-ts}dt = \sum_{k=0}^\infty f_{2k}\frac{\partial^{2k}L\varphi(s)}{\partial s^{2k}} - \sum_{k=0}^\infty f_{2k+1}\frac{\partial^{2k+1}L\varphi(s)}{\partial s^{2k+1}}, s>0.$$

Where $L(f)$ is the Laplace transform of the function $f$ .

**Proof.** Let $f(x)=\sum_{n=0}^\infty f_n x^n$ , is the series expansion of $f$ then

$$\int_0^\infty f(t)\varphi(t)e^{-ts}dt = \int_0^\infty\left(\sum_{k=0}^\infty f_{2k}t^{2k}\right)\varphi(t)e^{-ts}dt + \int_0^\infty\left(\sum_{k=0}^\infty f_{2k+1}t^{2k+1}\right)\varphi(t)e^{-ts}dt =$$

$$=\sum_{k=0}^\infty f_{2k}\int_0^\infty t^{2k}\varphi(t)e^{-ts}dt - \sum_{k=0}^\infty f_{2k+1}\int_0^\infty t^{2k+1}\varphi(t)e^{-ts}dt =$$

$$\sum_{k=0}^\infty f_{2k}\frac{\partial^{2k}}{\partial s^{2k}}\int_0^\infty \varphi(t)e^{-ts}dt - \sum_{k=0}^\infty f_{2k+1}\frac{\partial^{2k+1}}{\partial s^{2k+1}}\int_0^\infty \varphi(t)e^{-ts}dt =$$

$$=\sum_{k=0}^\infty f_{2k}\frac{\partial^{2k}}{\partial s^{2k}}L(\varphi)(s) - \sum_{k=0}^\infty f_{2k+1}\frac{\partial^{2k+1}}{\partial s^{2k+1}}L(\varphi)(s). \ \square$$

**Proposition 11. (Generalization of Proposition 1)**

Let $f(x)$ be as in Lemma, then there exists a unique function

$\psi_a(x) = L^{-1}\left(\log(\vartheta_4(it/2, e^{-a\pi}))\right)(x)$, such that:

$$\int_0^\infty f(t)\psi_a(t)e^{-ts}dt = \sum_{n=0}^\infty (-1)^n f_n \cdot \frac{\partial^n}{\partial s^n}\log(\vartheta_4(is/2, e^{-\pi a})) =$$

$$= f(0)\log\left(\prod_{n=1}^\infty\left(1-e^{-2n\pi a}\right)\right) - \sum_{n\in\mathbb{Z}-\{0\}}\frac{f(n)e^{-ns}}{2n\sinh(\pi a n)} \quad : (12)$$

$$\int_0^\infty f(e^{-t})\psi_a(t)e^{-ts}dt = \sum_{n=0}^\infty f_n \log(\vartheta_4(i(s+n)/2, e^{-\pi a})) =$$

$$f(1)\log\left(\prod_{n=1}^\infty\left(1-e^{-2n\pi a}\right)\right) - \sum_{n\in\mathbb{Z}-\{0\}}\frac{f(e^{-n})e^{-ns}}{2n\sinh(\pi a n)} : (13)$$

$a$ is a positive parameter.

**Proof.**

Let $f_n = \dfrac{f^{(n)}(0)}{n!}$ and $f(x)=\sum_{n=0}^\infty f_n x^n$ . We have:

$$\sum_{n=1}^\infty \frac{\cosh(tn)}{n\sinh(\pi a n)} = \log(P_0) - \log(\vartheta_4(it/2, e^{-a\pi})), \text{ when } |t| < a\pi$$

Differentiating with respect to $t$ we get:

$$\sum_{n=1}^{\infty}\frac{\cosh(tn)n^{2k}}{n\sinh(\pi an)}=-\frac{\partial^{2k}}{\partial t^{2k}}\log(\vartheta_4(it/2,e^{-a\pi})),$$

$$\sum_{n=1}^{\infty}\frac{\sinh(tn)n^{2k+1}}{n\sinh(\pi an)}=-\frac{\partial^{2k+1}}{\partial t^{2k+1}}\log(\vartheta_4(it/2,e^{-a\pi})),$$

when $k=0,1,2,\ldots$.

Also

$$\sum_{n=1}^{\infty}\frac{\cosh(tn)f_{2k}n^{2k}}{n\sinh(\pi an)}=-f_{2k}\frac{\partial^{2k}}{\partial t^{2k}}\log(\vartheta_4(it/2,e^{-a\pi}))$$

$$\sum_{n=1}^{\infty}\frac{\cosh(tn)f_{2k+1}(0)n^{2k+1}}{n\sinh(\pi an)}=-f_{2k+1}\frac{\partial^{2k+1}}{\partial t^{2k+1}}\log(\vartheta_4(it/2,e^{-a\pi}))$$

If we sum the above relations and use the lemma we have

$$\int_0^{\infty}f(t)\varphi(t)e^{-ts}dt=c-\sum_{n=1}^{\infty}\frac{\cosh(tn)}{n\sinh(\pi an)}f_e(n)+\sum_{n=1}^{\infty}\frac{\sinh(tn)}{n\sinh(\pi an)}f_o(n)$$

where $f_e(x)=\frac{f(x)+f(-x)}{2}$ and $f_o(x)=\frac{f(x)-f(-x)}{2}$ after some simplifications we complete the proof of the second equality of (12)

For the proof of the first equality of (12) observe that

$$\int_0^{\infty}(-t)^n\psi_a(t)e^{-ts}dt=\frac{\partial^n}{\partial s^n}\left(\log(\vartheta_4\left(is/2,e^{-a\pi}\right)\right)$$

and use $f(x)=\sum_{n=0}^{\infty}f_nx^n$ . □

In the same way we have

**Proposition 12.**
Let $f(x)$ be as in Lemma, then there exists a unique function $\rho_a(x)=L^{-1}\left(\log(\vartheta_2(t,e^{-1/a}))\right)(x)$, such that:

$$\int_0^{\infty}f(t)\rho_a(t)e^{-ts}dt=\sum_{n=0}^{\infty}(-1)^nf_n\cdot\frac{\partial^n}{\partial s^n}\log(\vartheta_2(s,e^{-1/a}))=$$

$$=2a-2af(0)s+a\pi\sum_{n\in\mathbb{Z}-\{0\}}\frac{f(2\pi na)e^{-2\pi nsa}}{\sinh(\pi^2an)}\quad :(14)$$

If all sides of (14) are convergent and $a$ is a positive parameter.